 \documentclass[preprint,12pt]{elsarticle}
\usepackage[centertags]{amsmath}
\usepackage{amsfonts}
\usepackage{amssymb}
\usepackage{amsthm}
\usepackage{newlfont}
\usepackage{graphicx}

\newtheorem{theorem}{Theorem}[section]
\newtheorem{lemma}[theorem]{Lemma}
\newtheorem{corollary}[theorem]{Corollary}

\theoremstyle{definition}

\usepackage{amssymb}

 \journal{}

\begin{document}

 \begin{frontmatter}



  \fntext[lable3]{Corresponding
author. \textit{E-mail
addresses: \small{m$_-$nabiei@sbu.ac.ir}}}

 \title{Newton type method for nonlinear Fredholm
integral equations}

 \author{\small{Mona Nabiei,$^{1}$ Sohrab Ali Yousefi.$^{2}$ 
}}
 \address{$^{1,2}$Department of Mathematics, Shahid Beheshti University, G. C. P.O. Box 19839, Tehran, IRAN
}

 \begin{abstract}
 This paper present a numerical method for solving nonlinear Fredholm integral equations.
 The method is based upon Newton type approximations. Illustrative examples are included to demonstrate the validity and applicability of the technique.
 \end{abstract}


 \end{frontmatter}


\section{Introduction and some preliminaries}

 The solution of the nonlinear Fredholm integral equation has been of considerable concern. This equation arises in the
theory of parabolic boundary value problems, the mathematical
modeling of the spatio-temporal development of an epidemic and
various physical and biological problems. A through discussion of
the formulation of these models and the scientific structure  are
given in \cite{7,22,23} and the references therein.

The nonlinear Fredholm integral equation is given in \cite{23} as
\begin{equation}\label{1}
h(x)=f(x)+\lambda \int_{0}^{1}G(x,t,h(t))dt,
\end{equation}

 where  $h(x)$, is an unknown function,  the functions $f(x)$ and
$G(x,t,h(t))$ are analytic on $[0,1]$. The Existence and
uniqueness results for Eq.(\ref{1}) may be found in \cite{7,11,14}.
However, few numerical methods for Eq.(\ref{1}) are known in the
literature \cite{23}. For the linear case, the time collocation
method was introduced in \cite{17} and projection method was
presented in \cite{12,13}. In \cite{2},  the results of \cite{17}
have been extended to nonlinear Volterra-Hammerstein integral
equations. In \cite{15,23}, a technique based on the Adomian
decomposition method was used for the solution of  Eq.(\ref{1}).

The next contexts from \cite{19} will be used later on.
Recall that the lower sub differential of a function $f:X \rightarrow [ - \infty , +
\infty ]$ on a Banach space $X$ at $x\in D(f)$, is given by
$$ \partial f(x)=\{ x^* \in X^* : ~~ \forall \nu \in X ~~ < x^* , \nu >  \leq f'(x,\nu) \},$$
where $X^*$ is the dual space of $X$ and $f'$ is the lower sub
derivative of $f$ at  $x$, given by
$$ f'(x,\nu)= \liminf_{\begin{array}{c}
 t\rightarrow 0\\
 u\rightarrow \nu
\end{array}} \frac{f(x+tu)-f(x)}{t},~~~~~~~~~\nu \in X. $$
And $\partial f$ is given by
$$\partial f =\{(x,x^*):~~~ x\in D(f),~~x^*\in \partial f(x)\}. $$
A Banach space $X$ is said to be a reliable space, in short an
R-space, if for any lower semi continuous function $f:X
\rightarrow [ - \infty , + \infty ]$, for any convex Lipschitzian
function $g$ on $X$, for any $x\in D(f)$ at which $f+g$ attains its
infimum, and for any $\varepsilon >0$, one has $0\in
\partial f(u)+ \partial g(v)+ \varepsilon B^*,$ for some
$u,v\in B_{\varepsilon}(x)$ such that $|f(u)-f(x)|<\varepsilon$.

We say that a function $f:X \rightarrow ( - \infty , + \infty )$
on a Banach space $X$ is lower $T^1$ or of class $LT^1$ if it is
lipschitzian, $G\hat{a}teaux$  differentiable, and if its
derivative $f':X\rightarrow X^*$ is continuous when $X^*$ is
endowed with the $weak^*$ topology.

A Banach space $X$ is said to be $LT^1$-bumpable if there exists
a nonull function of class $LT^1$ on $X$ with bounded support.
Any Banach space whose norm is $G\hat{a}teaux$ differentiable
off $0$ is $LT^1$-bumpble. The fact that any $LT^1$-bumpable
Banach space is reliable is used in the following.

\begin{theorem}\label{r}
Let $f:X\rightarrow [-\infty,+\infty ]$ be a lower semi continuous on an R-space. If there exists a constant $k\geq 0$ such that for any $(x,x^*)\in \partial f$ one has $\|x^*\|\leq k$, then $f$ is Lipschitzian with rate $k$ on its domain $D(f)$.
\end{theorem}

\section{Newton type method}

Consider the nonlinear Fredholm integral equation is given in
\cite{23} as
$$h(x)=f(x)+\lambda \int_{0}^{1}G(x,t,h(t))dt.$$

Let $f$ be continuous on $[0,1]$, and continuous function $G$ have continuous, bounded
derivative on it's third component.
 
We define operator
$F:C[0,1]\rightarrow C[0,1]$ by
$$ F(h)(x)=h(x)-f(x)-\lambda \int_{0}^{1}G(x,t,h(t))dt,~~~~~~~~ h\in C[0,1],~x\in [0,1],$$
and for each $u\in C[0,1]$, define operator
$T_{F,u}:C[0,1]\rightarrow C[0,1]$ as:
$$
T_{F,u}(h)(x)=\lim_{\varepsilon\rightarrow 0}\frac{F(h+\varepsilon u)(x)-F(h)(x)}{\varepsilon}~~~~~~~~h\in C[0,1],~x\in [0,1].$$
So we have
$$T_{F,u}(h)(x)=u(x)-\lambda \int_{0}^{1}\frac{\partial G(x,t,h(t))}{\partial h}u(t)dt.$$

\textbf{Note.} The Operator $T_{F,u}$ is continuous with respect
to $h$ and $u$. Put
$$H_u(h)(x)=h(x)-\frac{F(h)(x)}{T_{F,u}(h)(x)}, ~~~~h,u\in C[0,1], x\in [0,1].$$ 
We say that $F$ is $u$-smooth, if $T_{F,u}(h)(x)$ is nonzero for each $h\in C[0,1]$ and $x\in [0,1]$. If $F$ is $u$-smooth then the operator $H_u:C[0,1]\rightarrow C[0,1]$ is well-define.

\textbf{Remark.} Let $F$ be $1$-smooth. If there exists $p\in C[0,1]$ such that $F(p)=0$,
then for each $u\in C[0,1]$ and $x\in [0,1]$, we have

\begin{eqnarray}
T_{H_1,u}(p)(x) &=& \lim_{\varepsilon\rightarrow 0}\frac{H_1(p+\varepsilon u)(x)-H_1(p)(x)}{\varepsilon}\nonumber\\
~&=&\lim_{\varepsilon\rightarrow 0}u(x)- \frac{F(p+\varepsilon u)(x)}{\varepsilon}\frac{1}{T_{F,1}(p+\varepsilon u)(x)}\nonumber\\
~&=&u(x)-T_{F,u}(p)(x)\frac{1}{T_{F,1}(p)(x)}.\nonumber
\end{eqnarray}

So, if $F$ is $1$-smooth and $F(p)=0$, then $T_{H_1,1}(p)=0$. Hence the hypothesis of
continuity, implies, there exist neighborhoods of $p$ and $1$,
respectively denoted by $B_r(p)$ and $B_\delta (1)$, for suitable
$r$ and $\delta$, such that
$$\sup_{x\in [0,1]}|T_{H_1,u}(h)(x)|\leq \frac{1}{2},$$
for each $ h\in B_r(p)$ and $u\in B_\delta (1)$.

Now for each $x\in [0,1]$, define $\varphi (x)$ from $C[0,1]$ to real line, with
$$\varphi (x)(h)=H_1(h)(x),$$
where $h$ is in $C[0,1]$, and consider the next Lemma.

\begin{lemma}
Let $F$ be $1$-smooth, $x\in [0,1]$, $h \in B_r(p)\subseteq C[0,1]$ and $h^*\in (C[0,1])^*$ are fix. If for each $g\in C[0,1]$ there exist a
neighborhood $B_\gamma (g)$ and $\varepsilon >0$, such that
$$h^*(g)\leq \frac{\varphi(x)(h+\varepsilon u)-\varphi (x)(h)}{\varepsilon},$$
for each $u \in B_\gamma (g)$, then $||h^*||\leq \frac{1}{2}$.
\end{lemma}
\textbf{proof.} By hypothesis, $h^*(g)\leq \frac{H_1(h+\varepsilon u)(x)-H_1(h)(x) }{\varepsilon}$, implies by attending
$\varepsilon$ to zero, $h^*(g)\leq T_{H_1,u}(h)(x)$, for each
$u\in B_\gamma(g)$. Now let $g\in C[0,1]$ and $||g-1||<
\frac{\delta}{2}$, there exists $\gamma_1>0$ such that for each $u\in B_{\gamma_1}(g)$,  
$h^*(g)\leq T_{H_1,u}(h)(x)$. Put
$\gamma=\min\{\gamma_1,\frac{\delta}{2}\}$, hence for each $u\in
B_\gamma (g)$, we have
\[ \|u-1\|\leq \|u-g\|+\|g-1\|<\frac{\delta}{2}+\frac{\delta}{2}=\delta,\]
and before Remark implies that
$h^*(g)<T_{H_1,u}(h)(x)\leq\frac{1}{2}$. If $||g||=1$ and
$||g-1||<\delta$ then $h^*(g)\leq\frac{1}{2}$. So we assume that
$\delta <||g-1||$, then the linearity of $h^*$, by using of the
convexity of closed unit ball in $C[0,1]$, implies that there
exist $\alpha\in (0,1)$ and $g'\in C[0,1]$, with
$||g'-1||<\delta$, $h^*(g')\leq h^*(1)$, and
$g=\frac{1}{\alpha}g'+(1-\frac{1}{\alpha})1$. Hence
$$h^*(g)=\frac{1}{\alpha}h^*(g')+(1-\frac{1}{\alpha})h^*(1)\leq \frac{1}{\alpha}h^*(1)+(1-\frac{1}{\alpha})h^*(1)=h^*(1).$$
But $1\in B_\delta (1)$ and $h^*(g)\leq \frac{1}{2}$. Using this
fact, by linearity of $h^*$ implies that there exists $g\in
C[0,1]$, with $||g||=1$ and $|h^*(g)|\leq\frac{1}{2}$, so
$||h^*||\leq \frac{1}{2}$ and our proof is completed.
$\blacksquare$

\textbf{Note.} By considering theorem 5.4 of \cite{6}, we know
that the $C[0,w_0]$, it can be endowed with Frechet differentiable
norm and hence it can be a $LT^1$-bumpable space. So it can be a
$R$-space. Because each $LT^1$-bumpable space is
reliable\cite{18,19}.

\begin{theorem}
If $x\in [0,1]$, $(h,h^*)\in \partial\varphi (x)$ and $||h-p||<r$,
then $||h^*||\leq\frac{1}{2}$.
\end{theorem}
\textbf{proof.} Since $(h,h^*)\in \partial\varphi (x)$, for each
$g\in C[0,1]$, we have $h^*(g)\leq \varphi'(x)(h,g)$ and before
Lemma implies that $||h^*||\leq \frac{1}{2}$. $\blacksquare$

\begin{corollary}
The function $\varphi (x)$ is Lipschitzian on $B_r(p)$ with
constant $\frac{1}{2}$, for each $x\in [0,1]$.
\end{corollary}
\textbf{proof.} Put $D\varphi (x)=B_r(p)$ as a subset of a
reliable space. So theorem\ref{r} immediately
implies the conclusion. $\blacksquare$

Now we are already to prove that the sequence
$$u_{n+1}(x)=u_n(x)-\frac{F(u_n)(x)}{T_{F,1}(u_n)(x)}~~~~~~x\in [0,1],~n=0,1,2,...$$
 converges to the solution of Eq. (1).

\begin{theorem}
Let $F$ be $1$-smooth, $p\in C[0,1]$ and $F(p)=0$, then there exists a neighborhood
of $p$ in $C[0,1]$ such that the above sequence $\{u_n\}$ is
convergence to $p$, with each arbitrary started point $u_0$  of
this neighborhood.
\end{theorem}
\textbf{proof.} Consider the neighborhood, which exposed to the
before Remark. Now for each $u_0$ of this neighborhood, we have
$||u_0-p||<r$, hence
$$
\begin{array}{ccc}
||u_n-p||&\leq & \sup_{x\in [0,1]}|H_1(u_{n-1})(x)-H_1(p)(x)|\\
~&=& \sup_{x\in [0,1]}|\varphi (x)(u_{n-1})-\varphi (x)(p)|~~~\\
~&\leq & \frac{1}{2}||u_{n-1}-p||~~~~~~~~~~~~~~~~~~~~~~~~~~~~~\\
~&\leq&...\leq{\frac{1}{2}}^n||u_0-p||<r{\frac{1}{2}}^n~~~~~~~~~~~~~~~~\\

\end{array}
$$
and the proof is completed. $\blacksquare$

\section{Illustrative Examples}

We applied the method  presented in this paper and solved some
test problems.

\subsection{Example 1}

Consider the nonlinear Fredholm integral equation

\begin{equation} \label{11}
u(x)=x^2-\frac{1}{8}\cos(1)+\frac{1}{8}-\frac{1}{4}\int_0^1 t \sin(u(t))dt,
\end{equation}

which has the exact solution $u(x)=x^2.$\\
Note that
\begin{eqnarray}
F(h)(x) &=& h(x)-x^2+\frac{1}{8}\cos(1)-\frac{1}{8}+\frac{1}{4}\int_0^1 t \sin(h(t))dt,\nonumber\\
T_{F,1}(h)(x) &=& 1+\frac{1}{4}\int_0^1t \cos(h(t))dt\neq 0.\nonumber
\end{eqnarray}

We applied the method  presented in this paper and solved Eq.(\ref{11})

The corresponding Newton sequence is
$$u_{n+1}(x)=u_n(x)-\frac{u_n(x)-(x^2-\frac{1}{8}\cos(1)+\frac{1}{8}-\frac{1}{4}\int_0^1 t\sin(u_n(t))dt)}{1+\frac{1}{4}\int_0^1 t\cos(u_n(t))dt}$$
and $u_0(x)=1.$

In figure 1 the exact solution and approximate solution $u_3(x)$
are plotted.

\begin{center}
\begin{figure}
\input{epsf}
\epsfxsize=7cm \epsfysize=7cm
\centerline{\epsffile{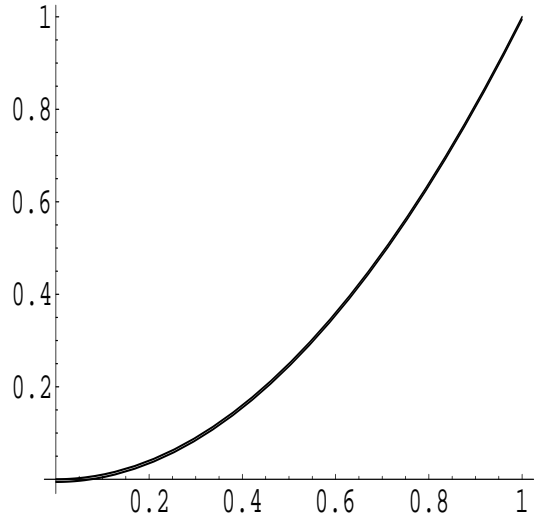}} \caption{Exact and
Approximate Solution of $u(x)$}
\end{figure}
\end{center}

\subsection{Example 2}

Consider the nonlinear Fredholm integral equation
\begin{equation} \label{22}
u(x)=e^x-\frac{1}{2}x(\cos(1)-\cos(e))+\frac{1}{2}\int_0^1xe^t \sin(u(t))dt, 
\end{equation}
which has the exact solution $u(x)=e^x$. \\
Note that
\begin{eqnarray}
F(h)(x) &=&h(x)-e^x+\frac{1}{2}x(\cos(1)-\cos(e))-\frac{1}{2}\int_0^1 x e^t \sin(h(t))dt,\nonumber\\
T_{F,1}(h)(x) &=&1-\frac{1}{2}\int_0^1 x e^t \cos(h(t))dt\neq 0.\nonumber
\end{eqnarray}

We applied the method  presented in this paper and solved Eq.(\ref{22})

The corresponding Newton sequence is
$$u_{n+1}(x)=u_n(x)-\frac{u_n(x)-(e^x-\frac{1}{2}x(\cos(1)-\cos(e))+\frac{1}{2}\int_0^1 x e^t\sin(u_n(t))dt}{1-\frac{1}{2}\int_0^1 x e^t \cos(u_n(t))dt}$$
and $u_0(x)=1.$

In figure 2 the exact and approximate solution $u_7(x)$ are
plotted.

\begin{center}
\begin{figure}
\input{epsf}
\epsfxsize=7cm \epsfysize=7cm
\centerline{\epsffile{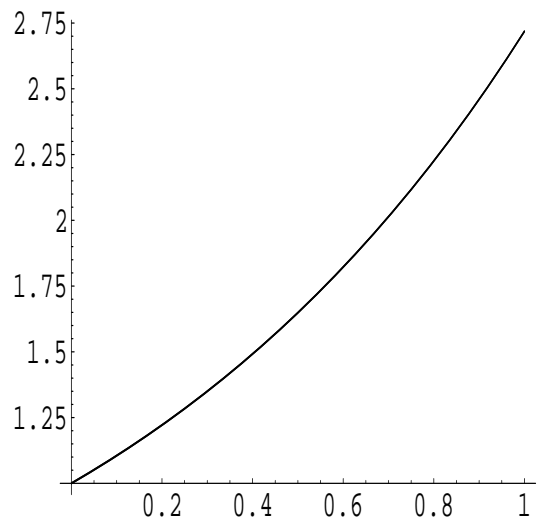}} \caption{Exact and
Approximate Solution of $u(x)$}
\end{figure}
\end{center}

\section{CONCLUSION}

The properties of
the Legendre wavelets together with the Gaussian integration
method are used to reduce the solution of the mixed Volterra-Fredholm integral equations to the solution of
algebraic equations. Illustrative examples are included to
demonstrate the validity and applicability of the technique. Moreover, only a small number of Legendre wavelets are needed to obtain a satisfactory result. The
given numerical examples support this claim.


\end{document}